\documentclass[11pt]{article}
\usepackage{graphicx,amsmath,amssymb,amsfonts,color}
\usepackage{pstcol}

\begin{document}

\def\nt{\noindent}

\centerline{\Large \bf On the $k$-edge magic graphs}

\bigskip
\baselineskip12truept \centerline{ Sin-Min Lee$^{a}$, Saeid
Alikhani$^{b,}$\footnote{\baselineskip12truept\it\small E-mail:
alikhani@yazduni.ac.ir}, Gee-Choon Lau$^c$ and William Kocay$^d$}
\bigskip
\baselineskip20truept \centerline{\it $^a$Department of Computer
Science} \vskip-8truept \centerline {\it San Jose State University
San Jose, California 95192 U.S.A.}

\medskip
 \vskip-8truept \centerline{\it $^{b}$Department of
Mathematics} \vskip-9truept \centerline{\it Yazd University,
89195-741, Yazd, Iran}
\medskip
 \vskip-9truept \centerline{\it $^c$Faculty of Computer and Mathematical Sciences}
\vskip-8truept\centerline{\it Universiti Teknologi MARA, Johor
85009 Segamat, Malaysia}

\medskip
 \vskip-9truept \centerline{\it $^d$Department of Computer Science
University of Manitoba} \vskip-8truept\centerline{\it Winnipeg,
Canada R3T 2N2}

\begin{abstract}
Let $G$ be a graph with vertex set $V$ and edge set $E$ such that
$|V| = p$ and $|E| = q$. For integers $k \geq 0$, define an edge
labeling $f: E\rightarrow \{k, k+1,\ldots, k+q-1\}$ and define the
vertex sum for a vertex $v$ as the sum of the labels of the edges
incident to $v$. If such an edge labeling induces a vertex
labeling in which every vertex has a constant vertex sum (mod p),
then $G$ is said to be $k$-edge magic ($k$-EM). In this paper, we
(i) show that all the maximal outerplanar graphs of order $p = 4,
5, 7$ are $k$-EM if and only if $k \equiv 2$ (mod p); (ii) obtain
all the maximal outerplanar graphs that are $k$-EM for $k = 3, 4$;
and (iii) characterize  all $(p, p-h)$-graph that are $k$-EM for
$h\geq 0$. We conjecture that all maximal outerplanar graphs of
prime order $p$ are $k$-EM if and only if $k\equiv 2$ (mod p).
\end{abstract}

\nt{\bf Keyword.} $k$-edge-magic, outer planar graph.

\end{document}